\documentclass[11pt]{article}
\renewcommand{\baselinestretch}{1.2}
\newcommand{\dated}{\mbox{} \hfill {\small [{\tt \today}]}} \usepackage{amsmath,amssymb,amscd}
\newcommand{\pf}[1]{\trivlist \item[\hskip\labelsep\it #1\ ]}
\newcommand{\varpf}[1]{\trivlist \item[\hskip\labelsep\sc #1:]}
\newcommand{\qedbox}{$\rlap{$\sqcap$}\sqcup$}
\newcommand{\qed}{\qquad \qedbox \endtrivlist}
\newcommand{\varqed}{\hfill \rule{0.6em}{0.6em} \endtrivlist}
\newenvironment{proof}{\pf{Proof}}{\qed}

\newenvironment{remark}{\pf{Remark}}{\endtrivlist}
\newenvironment{remarks}{\pf{Remarks} 
   \begin{enumerate}}{\end{enumerate} \endtrivlist}
\newenvironment{example}{\pf{Example}}{\endtrivlist}
\newenvironment{examples}{\pf{Examples} 
   \begin{enumerate}}{\end{enumerate} \endtrivlist}
\newenvironment{items}{
  \begin{enumerate} 
                    
  }{\end{enumerate}}
\newenvironment{alphitems}{
  \begin{enumerate} 
                    
  }{\end{enumerate}}
\newenvironment{keywords}{\noindent\small {\it Keywords\/}:}{\vskip 4pt}
\newenvironment{classification}{\noindent\small 2000 {\it Mathematics Subject
Classification\/}:}{\vskip 12pt}

\renewcommand{\implies}{\quad\Longrightarrow\quad}
\renewcommand{\iff}{\quad\Longleftrightarrow\quad}

\newcommand{\comps}{{\mathbb C}}
\newcommand{\reals}{{\mathbb R}}

\newcommand{\void}{\varnothing}
\newcommand{\tensor}{\otimes}

\newcommand{\Tensor}{\hat{\otimes}}

\newcommand{\sdprod}{\rtimes}
\newcommand{\cstar}{{C^\ast}}

\newcommand{\A}{{\mathfrak A}}
\newcommand{\B}{{\mathfrak B}}

\newcommand{\SL}{\operatorname{SL}}
\newcommand{\SIN}{\operatorname{SIN}}

\newcommand{\supp}{{\operatorname{supp}}}

\newtheorem{theorem}{Theorem}[section]
\newtheorem{lemma}[theorem]{Lemma}
\newtheorem{corollary}[theorem]{Corollary}
\newtheorem{proposition}[theorem]{Proposition}
\newtheorem{df}[theorem]{Definition}
\newenvironment{definition}{\begin{df} \rm}{\end{df}}

\newcommand{\WAP}{\mathcal{WAP}}
\newcommand{\AP}{\mathcal{AP}}
\newcommand{\LUC}{\mathcal{LUC}}
\newcommand{\SC}{\mathcal{SC}}
\newcommand{\SO}{\operatorname{SO}}
\title{A Connes-amenable, dual Banach algebra \\
need not have \\
a normal, virtual diagonal}
\author{\it Volker Runde}
\date{}
\begin{document}
\maketitle
\begin{abstract}
Let $G$ be a locally compact group, and let $\WAP(G)$ denote the space
of weakly almost periodic functions on $G$. We show that, if $G$ is a
$[\SIN]$-group, but not compact, then the dual Banach algebra $\WAP(G)^\ast$ does not have a normal, virtual diagonal.
Consequently, whenever $G$ is an amenable, non-compact $[\SIN]$-group, $\WAP(G)^\ast$ is an example of a Connes-amenable, dual Banach algebra without a normal, virtual diagonal.
On the other hand, there are amenable, non-compact, locally compact groups $G$ such that $\WAP(G)^\ast$ does have a normal, virtual diagonal.
\end{abstract}
\section*{Introduction}
In \cite{Joh1}, B.\ E.\ Johnson showed that a locally compact group $G$ is amenable if and only if its group algebra $L^1(G)$ has vanishing first order Hochschild cohomology with coefficients in dual
Banach $L^1(G)$-bimodules. Consequently, he called a Banach algebra satisfying this cohomological triviality condition amenable. Soon thereafter, Johnson gave a more intrinsic characterization of the amenable 
Banach algebras in terms of approximate and virtual diagonals (\cite{Joh2}). 
\par
For some classes of Banach algebras, amenability in the sense of \cite{Joh1} is too strong to allow for the development of a rich theory: it follows from work by S.\ Wassermann (\cite{Was1}), for example, that a von Neumann algebra 
is amenable if and only if it is subhomogeneous. This indicates that the definition of amenability should be modified when it comes to dealing with von Neumann algebras.
\par 
A variant of Johnson's definition that takes the dual space structure of a von Neumann algebra into account was introduced in \cite{JKR}, but is most commonly associated with A.\ Connes' paper \cite{Con1}. 
Following A.\ Ya.\ Helemski\u{\i} (\cite{Hel}), we shall refer to this variant of amenability as to Connes-amenability. As it turns out, Connes-amenability is equivalent to several other important properties of von
Neumann algebras, such as injectivity and semidiscreteness (\cite{BP}, \cite{Con1}, \cite{Con2}, \cite{EL}, \cite{Was2}; see \cite[Chapter 6]{LoA} for a self-contained exposition). Like the amenable Banach algebras,
the Connes-amenable von Neumann algebras allow for an intrinsic characterization in terms of diagonal type elements: a von Neumann algebra is Connes-amenable if and only if it has a normal, virtual diagonal (\cite{Eff}).
\par
The notions of Connes-amenability and normal, virtual diagonals make sense not only for von Neumann algebras, but for a larger class of Banach algebras called dual Banach algebras in \cite{Run1}. Examples of dual Banach algebras, besides von 
Neumann algebras, are --- among others --- the measure algebras $M(G)$
of locally compact groups $G$. As for von Neumann algebras,
amenability in the sense of \cite{Joh1} turns out to be too
restrictive a concept for measure algebras: the algebra $M(G)$ is
amenable if and only if $G$ is discrete and amenable (\cite{DGH}). In \cite{Run2},
however, the author showed that $M(G)$ is Connes-amenable if and only
if $G$ is amenable, and in \cite{Run3}, he also proved that
these conditions are also equivalent to $M(G)$ having a normal, virtual diagonal. 
\par
It is not hard to see that any dual Banach algebra with a normal, virtual diagonal must be Connes-amenable: as is observed in \cite{CG}, the argument from \cite{Eff} for von Neumann algebras carries over almost verbatim. The converse, however,
has been open so far (\cite[Problem 23]{LoA}).
\par
Besides $M(G)$ there are other dual Banach algebras associated with a locally compact group $G$. One of them is $\WAP(G)^\ast$, where $\WAP(G)$ denotes the weakly almost periodic functions on $G$. It is easy to see that
$\WAP(G)^\ast$ is Connes-amenable if $G$ is amenable; the converse is
also true, but not as straightforward (\cite{Run4}). If $G$ is compact, $\WAP(G)^\ast = M(G)$ has a normal, virtual diagonal, and in \cite{Run4}, the author made the --- 
as will become apparent: uneducated --- guess that $\WAP(G)^\ast$ has a normal, virtual diagonal if and only if $G$ is compact. 
\par
In the present paper, we shall confirm this conjecture for $[\SIN]$-groups. Consequently, whenever $G$ is an amenable $[\SIN]$-group that fails to be compact, the dual Banach algebra $\WAP(G)^\ast$ is Connes-amenable, but has no
normal, virtual diagonal. On the other hand, we shall see that $\WAP(G)^\ast$ does indeed have a normal, virtual diagonal if $G$ is amenable and minimally weakly almost periodic in the sense of \cite{Chou}. Since there are such
groups which fail to be compact, this shows that our conjecture from \cite{Run4} cannot be true in the generality stated there.
\subsection*{Acknowledgment}
I am grateful to Dona Strauss of Hull for bringing \cite{FS} to my attention and for providing Corollary \ref{idcor} and its proof.
\section{Connes-amenability and normal, virtual diagonals}
This section is preliminary in character: we briefly recall the definition of a dual Banach algebra along with the notions of Connes-amenability and of a normal, virtual diagonal.
\par
Given a Banach algebra $\A$ and a Banach $\A$-bimodule $E$, the dual space $E^\ast$ of $E$ becomes a Banach $\A$-bimodule in its own right via
\[
  \langle x, a \cdot \phi \rangle := \langle x \cdot a , \phi \rangle \quad\text{and}\quad  \langle x, \phi \cdot a \rangle := \langle a \cdot x  , \phi \rangle
  \qquad (a \in \A, \, \phi \in E^\ast, \, x \in E);
\]
in particular, the dual space $\A^\ast$ of $\A$ is a Banach $\A$-bimodule. Modules of this kind are referred to as {\it dual Banach modules\/}.
\par
The following definition was introduced in \cite{Run1}:
\begin{definition} \label{dualdef}
A Banach algebra $\A$ is called a {\it dual Banach algebra\/} if there is a closed submodule $\A_\ast$ of $\A^\ast$ such that $\A = (\A_\ast)^\ast$.
\end{definition}
\begin{remarks}
\item Equivalently, a Banach algebra $\A$ is dual if it is a dual Banach space such that multiplication is separately continuous in the $w^\ast$-topology.
\item In general, the predual space $\A_\ast$ in Definition \ref{dualdef} need not be unique, but will always be unambiguous from the context.
\end{remarks}
\begin{examples}
\item Every von Neumann algebra is a dual Banach algebra.
\item The measure algebra $M(G)$ of a locally compact group $G$ is a dual Banach algebra (with predual ${\cal C}_0(G))$.
\item If $E$ is a reflexive Banach space, then ${\cal B}(E)$ is a dual Banach algebra with predual $E \Tensor E^\ast$, where $\Tensor$ denotes the projective tensor product of Banach spaces.
\item The bidual of every Arens regular Banach algebra is a dual Banach algebra.
\end{examples}
\par
The following definition introduces a notion of amenability for dual Banach algebras that takes the dual space structure into account:
\begin{definition}
Let $\A$ be a dual Banach algebra.
\begin{alphitems}
\item A dual Banach $\A$-bimodule $E$ is called {\it normal\/} if the maps
\[
  \A \to E, \quad a \mapsto \left\{ \begin{array}{c} a \cdot x, \\ x \cdot a \end{array} \right.
\]
are $w^\ast$-continuous for each $x \in E$.
\item $\A$ is called {\it Connes-amenable\/} if every $w^\ast$-con\-ti\-nu\-ous derivation from $\A$ into a normal, dual Banach $\A$-bimodule is inner.
\end{alphitems}
\end{definition}
\begin{remarks}
\item For a von Neumann algebra, Connes-amenability is equivalent to a number of important properties, such as injectivity and semidiscreteness; see \cite[Chapter 6]{LoA} for a relatively self-contained account and for further references.
\item The measure algebra $M(G)$ of a locally compact group $G$ is Connes-amenable if and only if $G$ is amenable (\cite{Run2}).
\end{remarks}
\par
Let $\A$ be a Banach algebra. Then $\A \Tensor \A$ is a Banach $\A$-bimodule via
\[
  a \cdot (x \tensor y) := ax \tensor y \quad\text{and}\quad (x \tensor y) \cdot a := x \tensor ya \qquad (a,x,y \in \A),
\]
so that the multiplication map
\[
  \Delta \!: \A \Tensor \A \to \A, \quad a \tensor b \mapsto ab
\]
becomes a homomorphism of Banach $\A$-bimodules. Let $\A$ be a dual Banach algebra with predual $\A_\ast$, and let ${\cal B}_\sigma^2(\A,\comps)$ denote the bounded, bilinear functionals on $\A \times \A$
which are separately $w^\ast$-continuous, which form a closed submodule of $(\A \Tensor \A)^\ast$. Since $\Delta^\ast$ maps $\A_\ast$ into ${\cal B}_\sigma^2(\A,\comps)$, it follows that $\Delta^{\ast\ast}$ drops to an $\A$-bimodule homomorphism 
$\Delta_\sigma \!: {\cal B}_\sigma^2(\A,\comps)^\ast \to \A$. 
\par
We define:
\begin{definition} \label{normdiag}
A {\it normal, virtual diagonal\/} for a dual Banach algebra $\A$ is an element ${\mathrm M} \in {\cal B}_\sigma^2(\A,\comps)^\ast$ such that
\[
  a \cdot {\mathrm M} = {\mathrm M} \cdot a \quad\text{and}\qquad a \Delta_\sigma{\mathrm M} = a \qquad (a \in \A).
\]
\end{definition} 
\begin{remarks}
\item Every dual Banach algebra with a normal, virtual diagonal is Connes-amenable (\cite{CG}, but actually already \cite{Eff}).
\item A von Neumann algebra is Connes-amenable if and only if it has a normal, virtual diagonal (\cite{Eff}).
\item The same is true for the measure algebras of locally compact groups (\cite{Run2} and \cite{Run3}).
\end{remarks}
\par
As we shall see in this paper, there are Connes-amenable, dual Banach algebras which do not have a normal, virtual diagonal.
\section{The Banach algebra $\WAP(G)^\ast$ and the semigroup $G_\WAP$}
By a semitopological semigroup, we mean a semigroup $S$ equipped with a Hausdorff topology such that multiplication is separately continuous. If $S$ is locally compact, the measure space $M(S) \cong {\cal C}_0(S)^\ast$
can be turned into a Banach algebra via
\begin{equation} \label{conv}
  \langle f, \mu \ast \nu \rangle := \int_S \int_S f(st) \, d\mu(s) \, d\nu(t) \qquad (\mu , \nu \in M(S)).
\end{equation}
(Even though multiplication in $S$ need not be jointly continuous, the product integral on the right hand side of (\ref{conv}) does always exist and is independent of the order of integration; see \cite{Joh0}.) Note, that
even though $M(S)$ is a Banach algebra which is a dual Banach space, it need not be a dual Banach algebra in the sense of Definition \ref{dualdef}: this is due to the fact that ${\cal C}_0(S)$ need not be translation invariant. 
However, $M(S)$ {\it is\/} a dual Banach algebra if $S$ is compact or a group.
\par
Our first proposition is likely to be well known, but since we could not locate a reference, we include a proof:
\begin{proposition} \label{idprop}
Let $S$ be a locally compact, semitopological semigroup. Then the following are equivalent:
\begin{items}
\item $S$ has an identity.
\item $M(S)$ has an identity of norm one.
\end{items}
\end{proposition} 
\begin{proof}
(i) $\Longrightarrow$ (ii): If $S$ has an identity element, say $e$, then the point mass $\delta_e$ is an identity for $M(S)$ and trivially has norm one.
\par
(ii) $\Longrightarrow$ (i): Suppose that $M(S)$ has an identity element, say $\epsilon$, such that $\| \epsilon \| =1$; in particular,
\begin{equation} \label{ideq}
  \delta_s \ast \epsilon = \delta_s = \epsilon \ast \delta_s \qquad (s \in S)
\end{equation}
holds. It is straightforward that $\epsilon$ has to be an $\reals$-valued measure. Let $\epsilon^+, \epsilon^- \in M(S)$ be the Jordan decomposition of $\epsilon$, i.e.\ positive measures
such that $\epsilon = \epsilon^+ - \epsilon^-$ and $1 = \| \epsilon^+ \| + \| \epsilon^- \|$. Fix $s \in S$. Then $\delta_s \ast \epsilon^+$ and $\delta_s \ast \epsilon^-$ are positive measures such that
$\delta_s = \delta_s \ast \epsilon^+ - \delta_s \ast \epsilon^-$ and 
\[
  \| \delta_s \| = 1 = \| \epsilon^+ \| + \| \epsilon^- \| \geq \| \delta_s \ast \epsilon^+ \| + \| \delta_s \ast \epsilon^- \| \geq \| \delta_s \ast \epsilon \| = \| \delta_s \| = 1.
\]
The uniqueness of the Jordan decomposition of $\delta_s$ thus yields that $\delta_s \ast \epsilon^+ = \delta_s$. Analogously, one sees that $\epsilon^+ \ast \delta_s = \delta_s$.
Hence, (\ref{ideq}) still holds true if we replace $\epsilon$ by $\epsilon^+$.
\par
We shall see that the existence of a positive measure $\epsilon \in M(S)$ satisfying (\ref{ideq}) already necessitates $S$ to have an identity.
\par
Fix $s \in S$, and assume that there is $t \in \supp(\epsilon)$ such that $st \neq s$. We may choose a non-negative function $f \in {\cal C}_0(S)$ such that $f(s) = 0$ and $f(st) > 0$. We obtain that
\[
  0 < \int_S f(st) \, d\epsilon(t) = \langle f, \delta_s \ast \epsilon \rangle = \langle f, \delta_s \rangle = f(s) = 0,
\]
which is nonsense. Consequently, $st=s$ holds for all $t \in \supp(\epsilon)$; an analogous argument shows that $ts=s$ for all $t \in \supp(\epsilon)$. Since $s \in S$ was arbitrary, it follows that every element of
$\supp(\epsilon)$ is an identity for $S$.
\end{proof}
\par
Let $S$ be any semitopological semigroup, and let $f \!: S \to \comps$. For $s \in S$, we define the left translate $L_s f$ of $f$ by $s$ through
\[
  (L_s f)(t) := f(st) \qquad (t \in S).
\]
Let ${\cal C}_\mathrm{b}(S)$ denote the commutative $\cstar$-algebra of bounded, continuous functions on $S$.
\begin{definition}
Let $S$ be a semitopological semigroup. A bounded, continuous function $f \in {\cal C}_\mathrm{b}(S)$ is called {\it weakly almost periodic\/} if $\{ L_s f : s \in S \}$ is relatively compact in the weak topology on ${\cal C}_\mathrm{b}(S)$.
\end{definition}
\par
For any semitopological semigroup $S$, let
\[
  \WAP(S) := \{ f \in {\cal C}_\mathrm{b}(S) : \text{$f$ is weakly almost periodic} \}
\]
Our reference for almost periodic functions is mostly \cite{Bur}. It is easy to see that $\WAP(S)$ is a $\cstar$-subalgebra of ${\cal C}_\mathrm{b}(S)$ whose character space we denote by $S_\WAP$. It is clear, that $S_\WAP$ contains a canonical,
dense image of $S$. The multiplication of $S$ ``extends'' to $S_\WAP$, turning it into a compact, semitopological semigroup. For more on semigroup compactifications, see \cite{BJM}. The dual space $\WAP(S)^\ast$ can be identified with
$M(S_\WAP)$, and thus, in particular, becomes a dual Banach algebra. For an alternative definition of the multiplication on $\WAP(S)^\ast$, see \cite[(2.8) Corollary and (2.11) Proposition]{Pat}. 
\par
From now on, we shall only consider weakly almost periodic functions on locally compact groups. If $G$ is a locally compact group, ${\cal C}_0(G) \subset \WAP(G)$ holds and the canonical map from $G$ to $G_\WAP$ is
a homeomorphism onto its image (\cite[Theorem 3.6]{Bur}). It is straightforward that $G_\WAP \setminus G$ is a closed ideal of $G_\WAP$ and thus, in particular, is a compact, semitopological semigroup.
Let $\pi_0 \!: \WAP(G)^\ast \to M(G)$ be the restriction map from $\WAP(G)^\ast$ onto ${\cal C}_0(G)^\ast$. It is routinely checked that $\pi_0$ is a $w^\ast$-continuous algebra homomorphism. 
Consequently, ${\cal C}_0(G)^\perp = \ker \pi_0$ is a $w^\ast$-closed ideal in $\WAP(G)^\ast$ that can be identified, as a Banach algebra, with $M(G_\WAP \setminus G)$.
\par
As a corollary of Proposition \ref{idprop}, we obtain:
\begin{corollary} \label{idcor}
Let $G$ be a locally compact group. Then the following are equivalent:
\begin{items}
\item The ideal ${\cal C}_0(G)^\perp$ of $\WAP(G)^\ast$ has an identity.
\item The ideal $G_\WAP \setminus G$ of $G_\WAP$ has an identity.
\end{items}
\end{corollary}
\begin{proof}
All that needs to be shown is that, if ${\cal C}_0(G)^\perp$ has an identity element $\epsilon$, then $\| \epsilon \| = 1$ must hold. 
\par
Denote the bimodule module action of $\WAP(G)^\ast$ on $\WAP(G)$ by $\cdot$. Let $f \in \WAP(G)$ and observe that
\begin{eqnarray}
  | \langle f, \epsilon \rangle | & \leq & \sup \{ | \langle L_x f, \epsilon \rangle | : x \in G \}  \nonumber \\
  & = & \sup \{ | \langle f, \delta_x \ast \epsilon \rangle | : x \in G \} \nonumber \\
  & = & \sup \{ | \langle \epsilon \cdot f, \delta_x \rangle | : x \in G \} \nonumber \\
  & = & \| \epsilon \cdot f \|. \label{normone}
\end{eqnarray}
Since $\epsilon \in {\cal C}_0(G)^\perp$, the left hand side of (\ref{normone}) only depends on the equivalence class $\tilde{f}$ of $f$ in $\WAP(G) / {\cal C}_0(G) \cong {\cal C}(G_\WAP \setminus G)$.
We thus obtain:
\begin{eqnarray*}
  \left| \left\langle \tilde{f}, \epsilon \right\rangle \right| & \leq & \left\| \epsilon \cdot \tilde{f} \right\| \\
  & = & \sup \left\{ \left| \left\langle \epsilon \cdot \tilde{f}, \delta_s \right\rangle \right| : s \in G_\WAP \setminus G \right\} \\
  & = & \sup \left\{ \left| \left\langle \tilde{f}, \delta_s \ast \epsilon \right\rangle \right| : s \in G_\WAP \setminus G \right\} \\
  & = & \sup \left\{ \left| \left\langle \tilde{f}, \delta_s \right\rangle \right| : s \in G_\WAP \setminus G \right\} \\
  & = & \left\| \tilde{f} \right\|.
\end{eqnarray*}
Hence, $\| \epsilon \| \leq 1$ holds, which completes the proof.
\end{proof}
\par
In view of Corollary \ref{idcor}, we now turn to the question of whether, for a locally compact group $G$, the ideal $G_\WAP \setminus G$ can have an identity.
\par
For any locally compact group $G$, let $G_\LUC$ denote its $\LUC$-compactification (see \cite{BJM}). There is a canonical quotient map $\pi \!: G_\LUC \to G_\WAP$. An element $s \in G_\LUC$ is called a {\it point of unicity\/} if
$\pi^{-1}( \{ \pi(s ) \} ) = \{ s \}$.
\par
Recall that a locally compact group is called a $[\SIN]$-group if its
identity has a basis of neighborhoods invariant under conjugation; all
abelian, all compact, and all discrete groups are $[\SIN]$-groups.
\par
The following is (mostly) \cite[Theorem 1.4]{FS}:
\begin{theorem} \label{FSthm}
Let $G$ be a non-compact $[\SIN]$-group. Then $G_\LUC \setminus G$ contains a dense open subset $X$ consisting of points of unicity with the following properties:
\begin{items}
\item $\pi(X)$ is open in $G_\WAP \setminus G$;
\item $X$ is invariant under multiplication with elements from $G$;
\item $X$ has empty intersection with $(G_\LUC \setminus G)^2$.
\end{items}
\end{theorem}
\begin{remark}
Items (ii) and (iii) are not explicitly stated as a part of \cite[Theorem 1.4]{FS}, but follow from an inspection of the proof.
\end{remark}
\par
The following consequence of Theorem \ref{FSthm} was pointed out to me by Dona Strauss:
\begin{corollary} \label{Dona}
Let $G$ be a non-compact $[\SIN]$-group. Then the ideal $G_\WAP \setminus G$ does not have an identity.
\end{corollary}
\begin{proof}
Assume towards a contradiction that $G_\WAP \setminus G$ does have an identity, say $e$. Let $X$ be a set as specified in Theorem \ref{FSthm}. Since $\pi(X)$ is open in $G_\WAP \setminus G$, and since $Ge$ is dense in $G_\WAP \setminus G$,
it follows that $Ge \cap \pi(X) \neq \void$; from Theorem \ref{FSthm}(ii), we conclude that $e \in \pi(X)$. Let $p \in X$ be such that $\pi(p) = e$. Since $\pi(p^2) = e^2 = e$ and since $p$ is a point of unicity, it follows that
$p^2 = p$. This, however, contradicts Theorem \ref{FSthm}(iii).
\end{proof}
\section{Normal, virtual diagonals for $\WAP(G)^\ast$}
Given a locally compact group $G$, let $\SC(G_\WAP \times G_\WAP)$
denote the bounded, {\it separately\/} continuous functions on $G_\WAP
\times G_\WAP$. The space $\SC(G_\WAP \times G_\WAP)$ can be
canonically identified with ${\cal B}^2_\sigma(\WAP(G)^\ast,\comps)$ (\cite[Proposition 2.5]{Run2}). 
In terms of $\SC(G_\WAP \times G_\WAP)$, the bimodule action of
$\WAP(G)^\ast$ on ${\cal B}^2_\sigma(\WAP(G)^\ast,\comps)$ is given by
\[
  (\mu \cdot f)(s,t) := \int_{G_\WAP} f(s, tr) \, d\mu(r) \qquad (s,t \in G_\WAP)
\]
and
\[
  (f \cdot \mu)(s,t) := \int_{G_\WAP} f(rs, t) \, d\mu(r) \qquad (s,t \in G_\WAP)
\]
for $f \in \SC(G_\WAP \times G_\WAP)$ and $\mu \in \WAP(G)^\ast$ (this
is seen as in \cite[Proposition 3.1]{Run2}). 
\par
The verification of our first lemma in this section is routine:
\begin{lemma} \label{l1}
Let $G$ be a locally compact group, and let
\begin{equation} \label{iddef}
  I := \left\{ f \in \SC(G_\WAP \times G_\WAP) : \text{$f(s, \cdot) \in {\cal C}_0(G)$ for all $s \in G_\WAP$} \right\}.  
\end{equation}
Then:
\begin{items}
\item $\SC(G_\WAP \times G_\WAP)$, equipped with the supremum norm, is a commutative $\cstar$-algebra with identity;
\item $I$ is a closed ideal and a $\WAP(G)^\ast$-submodule of $\SC(G_\WAP \times G_\WAP)$.
\end{items}
\end{lemma}
\par
Let $\A$ be a $\cstar$-algebra, and let $I$ be a closed ideal of $\A$. As is well known, the second dual $\A^{\ast\ast}$ is a von Neumann algebra --- the enveloping von Neumann algebra of $\A$ --- 
containing $I^{\ast\ast}$ as a $w^{\ast\ast}$ closed ideal. The identity $P$ of $I^{\ast\ast}$ is a central projection in $\A^{\ast\ast}$ such that $I^{\ast\ast} = P \A^{\ast\ast}$.
\par
We make use of these facts in the case where $\A = \SC(G_\WAP \times G_\WAP)$ for a locally compact group $G$ and $I$ is as in (\ref{iddef}).
\begin{lemma} \label{l1b}
Let $G$ be a locally compact group, let $I$ be as in\/ {\rm (\ref{iddef})}, and let $P \in \SC(G_\WAP \times G_\WAP)^{\ast\ast}$ be the identity of $I^{\ast\ast}$. Then $P \cdot \delta_s = P$ holds for all $s \in G_\WAP$.
\end{lemma}
\begin{proof}
For convenience, set $\A := \SC(G_\WAP \times G_\WAP)$, and let $\Omega$ be
the character space of $\A$, so that $\A \cong {\cal C}(\Omega)$ via the
Gelfand transform. Through point evaluation, $\Omega$ contains a dense copy
$G_\WAP \times G_\WAP$. (Since functions in $\A$ are only {\it separately\/}
continuous on $G_\WAP \times G_\WAP$, the canonical map from $G_\WAP \times
G_\WAP$ into $\Omega$ need not be continuous.) Since $I$ is a closed ideal of
$\A$, there is an open subset $U$ of $\Omega$ such that $I \cong {\cal
  C}_0(U)$. Point evaluation maps $G_\WAP \times G$ onto a dense subset of
$U$.
\par
We claim that $U$ is dense in $\Omega$. Since $G_\WAP \times G_\WAP$ is dense
in $\Omega$, it is sufficient to show that each point in $G_\WAP \times
G_\WAP$ can be approximated by a net from $U$. Fix $(s,t) \in G_\WAP \times
G_\WAP$, and let $f \in \A$. Since $G$ is dense in $G_\WAP$, there is a net
$(x_\alpha)_\alpha$ in $G$ such that $x_\alpha \to t$; since $f(s, \cdot)$ is
continuous, we have $f(s,x_\alpha) \to f(s,t)$; and since $f \in \A$ is
arbitrary, this yields that $(s,x_\alpha) \to (s,t)$ in $\Omega$.
\par
Since $U$ is dense in $\Omega$, the ideal $I$ is essential in $\A$,
i.e.\ the only $f \in \A$ such that $f I = \{ 0 \}$ is $f = 0$. By the
universal property of the multiplier algebra, $\A$ thus canonically embeds
into ${\cal M}(I) \cong {\cal C}_{\mathrm b}(U)$, the multiplier algebra of $I$, which, in turn, can be
identified with the idealizer of $I$ in $I^{\ast\ast}$: $\{ F \in
I^{\ast\ast}: FI \subset I\}$ (for all this, see \cite{Ped}, for instance).
\par
All in all, we have a canonical, injective $^\ast$-homomorphism $\theta \!:
{\cal C}(\Omega) \to I^{\ast\ast}$ with $\theta({\cal C}(\Omega)) \subset {\cal M}(I)$, which is routinely seen
to satisfy 
\[
  \theta(f \cdot \delta_s) = \theta(f) \cdot \delta_s \qquad (f \in \A, \, s \in G_\WAP)
\]
(just test both sides against points in $G_\WAP \times G$). Since clearly
$\theta(1) = P$, we finally obtain that
\[ 
  P \cdot \delta_s = \theta(1) \cdot \delta_s = \theta(1 \cdot \delta_s) =
  \theta(1) = P \qquad (s \in G_\WAP).
\]
This completes the proof.
\end{proof}
\par
Since $\SC(G_\WAP \times G_\WAP)^{\ast\ast}$ is a von Neumann algebra for any locally compact group $G$, the canonical bimodule action of the commutative $\cstar$-algebra $\SC(G_\WAP \times G_\WAP)$ on its 
dual extends to a bimodule action of $\SC(G_\WAP \times G_\WAP)^{\ast\ast}$ on its predual $\SC(G_\WAP \times G_\WAP)^\ast$. We denote this module action by mere juxtaposition.
\begin{lemma} \label{l2}
Let $G$ be a locally compact group, let $I$ be defined as in\/ {\rm (\ref{iddef})}, and let $P$ denote the identity of $I^{\ast\ast}$. Then 
\[
  \delta_s \cdot (P \mathrm{N}) = P(\delta_s \cdot \mathrm{N}) \qquad (s \in G_\WAP))
\]
holds for all $\mathrm{N} \in {\cal B}_\sigma^2(\WAP(G)^\ast,\comps)^\ast$
\end{lemma}
\begin{proof}
First note that 
\begin{equation} \label{homeq}
  (fg) \cdot \delta_s = (f \cdot \delta_s)(g \cdot
  \delta_s) \qquad (s \in G_\WAP)
\end{equation}
for all $f, g \in \SC(G_\WAP \times G_\WAP)$. By separate $w^\ast$-continuity, it follows that (\ref{homeq}) holds as well for all $f,g
\in \SC(G_\WAP \times G_\WAP)^{\ast\ast}$.
\par
Fix $s \in G_\WAP$ and $\mathrm{N} \in {\cal
  B}_\sigma^2(\WAP(G)^\ast,\comps)^\ast$, and let $f \in \SC(G_\WAP
\times G_\WAP)$. We obtain:
\begin{eqnarray*}
  \langle f,\delta_s \cdot (P \mathrm{N}) \rangle & = & \langle f \cdot
  \delta_s, P \mathrm{N} \rangle \\
  & = & \langle (f \cdot \delta_s) P, \mathrm{N} \rangle \\
  & = & \langle (f \cdot \delta_s) (P \cdot \delta_s) , \mathrm{N}
  \rangle, \qquad \text{by Lemma \ref{l1b}}, \\  
  & = & \langle (fP) \cdot \delta_s, \mathrm{N} \rangle,
  \qquad\text{by (\ref{homeq})}, \\
  & = & \langle f P, \delta_s \cdot \mathrm{N} \rangle \\
  & = & \langle f, P(\delta_s \cdot \mathrm{N}) \rangle.
\end{eqnarray*}
This proves the claim.
\end{proof}
\begin{lemma} \label{l3}
Let $G$ be a locally compact group, such that $\WAP(G)^\ast$ has a normal, virtual diagonal, say $\mathrm{M}$, let $I$ be as in\/ {\rm (\ref{iddef})}, and let $P$ be the identity of $I^{\ast\ast}$. Define
\[
  \rho \!: \WAP(G)^\ast \to \WAP(G)^\ast, \quad \mu \mapsto
  \Delta_\sigma(P(\mathrm{M} \cdot \mu)).
\]
Then:
\begin{items}
\item ${\cal C}_0(G)^\perp$ is contained in $\ker \rho$;
\item $\pi_0 \circ \rho = \pi_0$ holds, where $\pi_0 \!: \WAP(G)^\ast \to M(G)$ is the canonical restriction map;
\item $\rho(\delta_s \ast \mu) = \delta_s \ast \rho(\mu)$ for all $s \in G_\WAP$ and $\mu \in \WAP(G)^\ast$.
\end{items}
\end{lemma}
\begin{proof}
(i): Let $\mu \in {\cal C}_0(G)^\perp$. It follows that $\mu \cdot f = 0$ for each $f \in I$ and, consequently, that $\mathrm{M} \cdot \mu \in I^\perp$. Let $F \in \SC(G_\WAP \times G_\WAP)^{\ast\ast}$, and note that $FP \in 
I^{\ast\ast}$. Since
\[
  \langle F, P(\mathrm{M} \cdot \mu) \rangle = \langle FP, \mathrm{M} \cdot \mu \rangle = 0,
\]
we conclude that $\mu \in \ker \rho$.
\par
(ii): Fix $\mu \in \WAP(G)^\ast$ and $f \in {\cal C}_0(G)$, and observe that
\begin{eqnarray*}
  \langle f, \rho(\mu) \rangle & = & \langle f, \Delta_\sigma(P(\mathrm{M} \cdot \mu)) \rangle \\
  & = & \langle \Delta^\ast f, P(\mathrm{M} \cdot \mu) \rangle \\
  & = & \langle (\Delta^\ast f)P, \mathrm{M} \cdot \mu \rangle \\
  & = & \langle \Delta^\ast f, \mathrm{M} \cdot \mu \rangle, \qquad\text{because $\Delta^\ast {\cal C}_0(G) \subset I$}, \\
  & = & \langle f, \Delta_\sigma(\mathrm{M} \cdot \mu) \rangle \\
  & = & \langle f , \mu \rangle.
\end{eqnarray*}
This proves the claim.
\par
(iii): Fix $\mu \in \WAP(G)^\ast$ and $s \in G_\WAP$. We obtain:
\begin{eqnarray*}
  \rho(\delta_s \ast \mu) & = & \Delta_\sigma(P((\mathrm{M} \cdot \delta_s) \cdot \mu)) \\
  & = & \Delta_\sigma(P(\delta_s \cdot (\mathrm{M} \cdot \mu ))), \qquad\text{by Definition \ref{normdiag}}, \\
  & = & \Delta_\sigma(\delta_s \cdot (P(\mathrm{M} \cdot \mu ))), \qquad\text{by Lemma \ref{l2}}, \\
  & = & \delta_s \ast \Delta_\sigma(P(\mathrm{M} \cdot \mu )) \\
  & = & \delta_s \ast \rho(\mu).
\end{eqnarray*}
This completes the proof.
\end{proof}
\par
We can now prove the main result of this section (and of the whole paper):
\begin{theorem} \label{mainthm}
Let $G$ be a non-compact $[\SIN]$-group. Then $\WAP(G)^\ast$ does not have a normal, virtual diagonal.
\end{theorem}
\begin{proof}
Assume towards a contradiction that $\WAP(G)^\ast$ has a normal, virtual diagonal. Let $\rho \!: M(G) \to \WAP(G)^\ast$ be as in Lemma \ref{l3}, and define $\epsilon := \delta_e - \rho(\delta_e)$, where $e$ is the identity of $G$. By Lemma \ref{l3}(ii), it is clear that
$\epsilon \in {\cal C}_0(G)^\perp$. Moreover, we have for $s \in G_\WAP \setminus G$ that
\begin{eqnarray*}
  \delta_s \ast \epsilon & = & \delta_s - \delta_s \ast \rho(\delta_e) \\
  & = & \delta_s - \rho(\delta_s \ast \delta_e), \qquad\text{by Lemma \ref{l3}(iii)}, \\
  & = & \delta_s, \qquad \text{by Lemma \ref{l3}(i)}.
\end{eqnarray*}
Consequently, $\epsilon$ is a right identity for ${\cal C}_0(G)^\perp$. 
\par
For any $f \in \WAP(G)$, the function $\check{f} \!: G \to \comps$ defined by letting $\check{f}(x) := f(x^{-1})$ for $x \in G$ lies also in $\WAP(G)$ (\cite[Corollary 1.18]{Bur}). Setting
\[
  \langle f, \check{\mu} \rangle := 
  \left\langle \check{f}, \mu \right\rangle 
  \qquad (\mu \in \WAP(G)^\ast, \, f \in \WAP(G) ),
\]
defines an anti-automorphism $\WAP(G)^\ast \ni \mu \mapsto \check{\mu}$ of $\WAP(G)^\ast$, which leaves ${\cal C}_0(G)^\perp$ invariant. Hence, $\check{\epsilon}$ is a left identity for ${\cal C}_0(G)^\perp$,
so that ${\cal C}_0(G)^\perp$ has in fact an identity. This, however, is not possible by Corollaries \ref{idcor} and \ref{Dona}.
\end{proof}
\par
The following corollary confirms the guess made at the end of \cite{Run4} for $[\SIN]$-groups:
\begin{corollary} \label{maincor}
Consider the following statements about a $[\SIN]$-group $G$:
\begin{items}
\item $\WAP(G)^\ast$ has a normal, virtual diagonal.
\item $G$ is compact.
\item $G$ is amenable.
\item $\WAP(G)^\ast$ is Connes-amenable.
\end{items}
Then
\[
  \text{\rm (i)} \iff \text{\rm (ii)} \implies \text{\rm (iii)} \iff \text{\rm (iv)}. 
\]
\end{corollary}
\begin{proof}
(i) $\Longrightarrow$ (ii) follows immediately from Theorem \ref{mainthm}, and the converse is shown in \cite{Run3} (for compact $G$, we have $\WAP(G)^\ast = M(G)$).
\par
(ii) $\Longrightarrow$ (iii) is well known (see \cite{Pat} or \cite[Chapter 1]{LoA}).
\par
(iii) $\Longleftrightarrow$ (iv) is \cite[Proposition 4.11]{Run4}.
\end{proof}
\par
Consequently, $\WAP(G)^\ast$ is a Connes-amenable dual Banach algebra
{\it without\/} a normal, virtual diagonal whenever $G$ is an
amenable, but not compact $[\SIN]$-group: this includes all
non-compact, abelian, locally compact groups as well as all infinite,
discrete, amenable groups.
\par
As a consequence of Corollary \ref{maincor}, we also obtain a
characterization of those locally compact groups $G$ for which
$\WAP(G)^\ast$ is amenable (in the sense of \cite{Joh1}):
\begin{corollary}
The following are equivalent for a locally compact group $G$:
\begin{items}
\item $\WAP(G)^\ast$ is amenable;
\item $G$ is finite.
\end{items}
\end{corollary}
\begin{proof}
Of course, only (i) $\Longrightarrow$ (ii) needs proof.
\par
If $\WAP(G)^\ast$ is amenable, so is its quotient $M(G)$. By
\cite{DGH}, this means that $G$ must be discrete. In particular, $G$
is a $[\SIN]$-group. Since $\WAP(G)^\ast$ is amenable, it must have a
normal, virtual diagonal, so that the discrete group $G$ is also compact by Corollary \ref{maincor}.
\end{proof}
\section{Minimally weakly almost periodic groups}
In view of Corollary \ref{maincor}, the conjecture (made in \cite{Run4}) is tempting that $\WAP(G)^\ast$ has a normal, virtual diagonal only if $G$ is compact. As we shall see in this final section, this is wrong.
\par
Recall that a continuous, bounded function $f$ on a locally compact group $G$ is called {\it almost periodic\/} if $\{ L_x f : x \in G \}$ is relatively compact in the norm topology of ${\cal C}_{\mathrm{b}}(G)$. Let
\[
  \AP(G) := \{ f \in {\cal C}_{\mathrm{b}}(G) : \text{$f$ is almost periodic} \}.
\]
Like $\WAP(G)$, the space $\AP(G)$ is a commutative $\cstar$-algebra. Its character space, denoted by $G_\AP$, is a compact group that contains a dense, but generally not
homeomorphic image of $G$ in a canonical manner. For more information, see \cite{Bur} or \cite{BJM}, for example.
\par
The following definition is from \cite{Chou}:
\begin{definition} \label{minwap}
Let $G$ be a locally compact group. Then $G$ is called {\it minimally weakly almost periodic\/} if $\WAP(G) = \AP(G) + {\cal C}_0(G)$.
\end{definition}
\begin{remarks}
\item Every compact group is trivially minimally weakly almost periodic. If $G$ is not compact, but weakly almost periodic, then the sum in Definition \ref{minwap} is a direct one.
\item The motion group $\reals^N \sdprod \SO(N)$ is minimally weakly almost periodic (and amenable) as is $\SL(2,\reals)$, which is not amenable (see \cite{Chou}).
\item If $G$ is minimally weakly almost periodic, the kernel of
  $G_\WAP$ (see \cite{Bur} for the definition) must equal $G_\WAP
  \setminus G$. Hence, by Corollary \ref{idcor}, a non-compact $[\SIN]$-group cannot be minimally weakly almost periodic.
(This follows also immediately from the main result of \cite{Chou2}).
\end{remarks}
\par
The verification of the following lemma is routine:
\begin{lemma} \label{diaglem}
Let $\A$ and $\B$ be dual Banach algebras each of which as a normal, virtual diagonal. Then $\A \oplus \B$ has a normal, virtual diagonal.
\end{lemma}
\par
It is now fairly straightforward to refute our ``conjecture'' from \cite{Run4}:
\begin{proposition} \label{dprop}
Let $G$ be a locally compact, minimally weakly almost periodic group. Then the following are equivalent:
\begin{items}
\item $G$ is amenable.
\item $\WAP(G)^\ast$ has a normal, virtual diagonal.
\end{items}
\end{proposition}
\begin{proof}
(i) $\Longrightarrow$ (ii): Without loss of generality, suppose that $G$ is not compact. Let $\pi_a \!: \WAP(G)^\ast \to \AP(G)^\ast$ and $\pi_0 \!: \WAP(G)^\ast \to M(G)$ be the respective restriction maps; they are
$w^\ast$-continuous algebra homomorphism. Since $\WAP(G) = \AP(G) \oplus {\cal C}_0(G)$, it follows that $\pi_a \oplus \pi_0 \!: \WAP(G)^\ast \to \AP(G)^\ast \oplus M(G)$ is a $w^\ast$-continuous isomorphism.
Since $G$ is amenable $M(G)$ has a normal, virtual diagonal by \cite{Run3}, and the same is true for $\AP(G)^\ast \cong M(G_\AP)$. From Lemma \ref{diaglem}, we conclude that $\WAP(G)^\ast$ has a normal, virtual diagonal.
\par
(ii) $\Longrightarrow$ (i) follows immediately from \cite[Proposition 4.11]{Run4}.
\end{proof}
\begin{example}
The motion group group $G := \reals^N \sdprod \SO(N)$ is minimally weakly almost periodic and amenable, so that $\WAP(G)^\ast$ has a normal, virtual diagonal even though $G$ fails to be compact.
\end{example}
\par
In view of Proposition \ref{dprop} and Corollary \ref{maincor}, the conjecture isn't farfetched that $\WAP(G)^\ast$ has a normal, virtual diagonal if and only if $G$ is amenable and minimally weakly almost periodic.
\dated
\vfill
\begin{tabbing} 
{\it Address\/}: \= Department of Mathematical and Statistical Sciences \\
                 \> University of Alberta \\
                 \> Edmonton, Alberta \\
                 \> Canada T6G 2G1 \\[\medskipamount]
{\it E-mail\/}:  \> {\tt vrunde@ualberta.ca}\\[\medskipamount]
{\it URL\/}:     \> {\tt http://www.math.ualberta.ca/$^\sim$runde/}
\end{tabbing}

\end{document}